\documentclass[a4paper,leqno,12pt]{amsart}

\usepackage{epsf,graphicx,epsfig}
\usepackage{amscd}
\usepackage{amsmath,latexsym,amssymb,amsthm}
\usepackage[nospace,noadjust]{cite}
\usepackage{textcomp}
\usepackage{setspace,cite}
\usepackage{lscape,fancyhdr,fancybox}
\usepackage[all,cmtip]{xy}
\renewcommand{\baselinestretch}{1.2}

\usepackage{graphicx}

\usepackage{color}
\usepackage{url}
\usepackage{enumerate}
\usepackage[mathscr]{euscript}

  \theoremstyle{plain}

\renewcommand{\thefootnote}{}
\newcommand{\bd}{\begin{document}}
	\newcommand{\ed}{\end{document}}
\newcommand{\bc}{\begin{center}}
	\newcommand{\ec}{\end{center}}
\newcommand{\vs}{\vspace}
\newcommand{\hs}{\hspace}
\newcommand{\bq}{\begin{quote}}
	\newcommand{\eq}{\end{quote}}
\newcommand{\mb}{\makebox}
\newcommand{\lt}{\left}
\newcommand{\rt}{\right}
\newcommand{\beqa}{\begin{eqnarray*}}
	\newcommand{\eeqa}{\end{eqnarray*}}
\newcommand{\beqn}{\begin{eqnarray}}
	\newcommand{\eeqn}{\end{eqnarray}}
\newcommand{\bbibl}{\begin{thebibliography}{99}}
	\newcommand{\ebibl}{\end{thebibliography}}
\newcommand{\ti}{\times}
\newcommand{\bit}{\begin{itemize}}
	\newcommand{\eit}{\end{itemize}}
\newcommand{\ben}{\begin{enumerate}}
	\newcommand{\een}{\end{enumerate}}
\newcommand{\lb}{\label}
\newcommand{\hf}{\hspace*{\fill}}
\newcommand{\vf}{\vspace*{\fill}}
\newcommand{\beq}{\begin{equation}}
	\newcommand{\eeq}{\end{equation}}
\newcommand{\ba}{\begin{array}}
	\newcommand{\ea}{\end{array}}
\newcommand{\del}{\partial}
\newcommand{\bm}[1]{\mb{\boldmath ${#1}$}}
\newcommand{\ot}{\otimes}
\newcommand{\nn}{\nonumber}
\newcommand{\R}{\mb{$I\!\!R$}}
\newcommand{\C}{{\mathcal C}}
\newcommand{\M}{{\mathcall M}}
\newcommand{\E}{{\mathcal E}}
\newcommand{\N}{{\mathcal N}}
\newcommand{\B}{{\mathcal B}}
\newcommand{\Y}{{\mathcal Y}}
\newcommand{\F}{{\mathcal F}}
\newcommand{\Rc}{{\mathcal R}}
\newcommand{\A}{{\mathcal A}}
\renewcommand{\P}{{\mathcal P}}
\renewcommand{\S}{{\mathcal S}}
\newcommand{\es}{\emptyset}
\newcommand{\ci}{\subseteq}
\newcommand{\cs}{\supseteq}
\renewcommand{\u}{\cup}
\renewcommand{\i}{\cap}
\newcommand{\bu}{\bigcup}
\newcommand{\bi}{\bigcap}
\newcommand{\la}{\leftarrow}
\newcommand{\ra}{\rightarrow}
\newcommand{\Ra}{\Rightarrow}
\newcommand{\Lra}{\Leftrightarrow}
\newcommand{\lgra}{\longrightarrow}
\newcommand{\Lgra}{\Longrightarrow}
\newcommand{\lglra}{\longleftrightarrow}
\newcommand{\Lglra}{\Longleftrightarrow}
\renewcommand{\a}{\alpha}
\renewcommand{\b}{\beta}
\newcommand{\g}{\gamma}
\newcommand{\G}{\Gamma}
\renewcommand{\d}{\delta}
\newcommand{\D}{\Delta}
\newcommand{\e}{\varepsilon}
\newcommand{\eps}{\epsilon}
\newcommand{\h}{\eta}
\renewcommand{\l}{\lambda}
\newcommand{\m}{\mu}
\newcommand{\n}{\nu}
\newcommand{\p}{\pi}
\newcommand{\s}{\sigma}
\newcommand{\Si}{\Sigma}
\newcommand{\ta}{\tau}
\newcommand{\ph}{\phi}
\newcommand{\Ph}{\Phi}
\renewcommand{\c}{\chi}
\newcommand{\om}{\omega}
\newcommand{\Om}{\Omega}
\newcommand{\tri}{\triangle}
\newcommand{\rec}[1]{\frac{1}{#1}}
\newcommand{\f}{\frac}
\newcommand{\sm}[2]{\sum_{#1}^{#2}}
\newcommand{\ld}{\ldots}
\newcommand{\ov}{\overline}
\newcommand{\ol}[1]{$\bar{\mb{#1}}$}
\newcommand{\un}{\underline}
\newcommand{\iy}{\infty}

\newcommand{\wt}{\widetilde}
\newcommand{\ds}{\displaystyle}
\newcommand{\bdm}{\begin{displaymath}}
	\newcommand{\edm}{\end{displaymath}}

\newcommand{\nin}{\not\in}
\newcommand{\bt}{\begin{tabular}}
	\newcommand{\et}{\end{tabular}}
\newcommand{\alter}[2]{\lt\{ \ba {ll}#1 \\ #2 \ea \rt.}
\newcommand{\alt}[4]{\lt\{ \ba{ll}#1 & \mb{if \,\,}#2 \\ #3 & \mb{if
		\,\,}#4 \ea \rt.}
\newcommand{\altn}[4]{\lt\{ \ba{rl}#1 & \mb{if \,\,}#2 \\ #3 & \mb{if
		\,\,}#4 \ea \rt.}
\newcommand{\alto}[6]{ \lt\{ \ba{ll}#1 & \mb{if \,\,}#2 \\ #3 & \mb{if
		\,\,} #4 \\ #5 & \mb{if \,\,}#6 \ea \rt.}
\newcommand{\altero}[5]{\mb{$\lt\{ \ba {ll}#1 & \mb{if \,\,}#2 \\ #3 &
		\mb{if \,\,} #4 \\ #5 & \mb{otherwise} \ea \rt.$}}

\newcounter{cnt1}
\newcounter{cnt2}
\newcounter{cnt3}
\newcommand{\blr}{\begin{list}{$($\roman{cnt1}$)$} {\usecounter{cnt1}
			\setlength{\topsep}{0pt} \setlength{\itemsep}{0pt}}}
	\newcommand{\bla}{\begin{list}{$($\alph{cnt2}$)$} {\usecounter{cnt2}
				\setlength{\topsep}{0pt} \setlength{\itemsep}{0pt}}}
		\newcommand{\bln}{\begin{list}{$($\arabic{cnt3}$)$} {\usecounter{cnt3}
					\setlength{\topsep}{0pt} \setlength{\itemsep}{0pt}}}
			\newcommand{\el}{\end{list}}
		\newcommand{\no}{\noindent}
		\newtheorem{Thm}{Theorem}[section]
		\newtheorem{Lem}[Thm]{Lemma}
		\newtheorem{Prop}[Thm]{Proposition}
		\newtheorem{Def}[Thm]{Definition}
		\newtheorem{Exm}[Thm]{Example}
		\newtheorem{Rem}[Thm]{Remark}
		\newtheorem{Cor}[Thm]{Corollory}
		\renewcommand{\baselinestretch}{1}
		\newcommand{\ilim}{\mathop{\varprojlim}\limits}
		\newcommand{\dlim}{\mathop{\varinjlim}\limits}
\newenvironment{mysubsection}[2][]
{\begin{subsec}\begin{upshape}\begin{bfseries}{#2.}
\end{bfseries}{#1}}
{\end{upshape}\end{subsec}}

\newenvironment{myeq}[1][]
{\stepcounter{thm}\begin{equation}\tag{\thethm}{#1}}
{\end{equation}}

\renewcommand{\Re}{\operatorname{Re}}
\renewcommand{\ker}{\operatorname{Ker}}
\newcommand{\Coker}{\operatorname{Coker}}
\newcommand{\colim}{\operatorname{colim}}
\newcommand{\diag}{\operatorname{diag}}
\newcommand{\Sing}{\operatorname{Sing}}
\newcommand{\Maps}{\operatorname{Map}}
\newcommand{\Hom}{\operatorname{Hom}}
\newcommand{\Aut}{\operatorname{Aut}}
\newcommand{\ob}{\operatorname{Ob}}
\newcommand{\Id}{\operatorname{Id}}
\newcommand{\pr}{\operatorname{pr}}
\newcommand{\rel}{\operatorname{rel}}
\newcommand{\Z}{\mathbb{Z}}

\title{}
\author{}
\date{}
\usepackage{amssymb}
\def\theequation {\arabic{equation}}
\renewcommand{\thefootnote}{}

\usepackage{hyperref}
\hypersetup{
	colorlinks,
	citecolor=blue,
	filecolor=black,
	linkcolor=blue,
	urlcolor=black
}

\begin{document}
			
			\title{A note on cohomology of Clifford algebras}
	
\author{Bikram Banerjee}
\email{pbikraman@rediffmail.com}
\address{Ranaghat College, Ranaghat, W.B. 741201, India.}

\author{Goutam Mukherjee}
\email{goutam.mukherjee@tcgcrest.org; gmukherjee.isi@gmail.com}
\address{Institute for Advancing Intelligence, TCG CREST, First Floor, Tower 1, Bengal Eco Intelligent Park (Techna), Block EM, Plot No 3, Sector V, Salt lake, Kolkata 700091, West Bengal, India.}

				
			\date{}
			\maketitle

			

			\thispagestyle{empty}
			\begin{abstract}
				In this article we construct a cochain complex of a complex Clifford algebra with coefficients in itself in a combinatorial fashion and we call the corresponding cohomology by {\it Clifford cohomology.} We show that {\it Clifford cohomology} controls the deformation of a complex Clifford algebra and can classify them up to Morita equivalence. We also study  Hochschild cohomology groups and formal deformations of the algebra of smooth sections of a complex Clifford algebra bundle over an even dimensional orientable Riemannian manifold \(M\) which admits a \(Spin^{c}\) structure.
			\end{abstract}

			\maketitle
			

			\vspace{0.5cm}

			\section{Introduction}
			
			Algebraic deformation theory of associative algebras was developed by M. Gerstenhaber \cite{gerstenhaber} \cite{gerstenhaber 2}. In algebraic deformation theory of associative algebras the Hochschild cohomogy group \(HH^{\ast}(A)\) of an associative algebra \(A\) with coefficients in \(A\) plays the key role. The second Hochschild cohomology group \(HH^{2}(A)\) of \(A\) has a one-to-one correspondence to the set of all equivalence classes of non isomorphic infinitesimal deformations of \(A\) while the triviality of \(HH^{2}(A)\) implies that any formal deformation of \(A\) is equivalent to null deformation. Again \(HH^{3}(A)= 0\) forces that any infinitesimal deformation of \(A\) can be extended to a formal deformation. Thus primarily the second and third Hochschild cohomology groups of an associative algebra \(A\) controls the deformations of \(A\). \newline
			
			Now let \(Ass\) denotes the category of associative algebras over some field \(k\) of characteristic zero, \(\mathcal{A}\) be a particular subcategory of \(Ass\) over  \(k\) having some extra structure. It is natural to ask whether the algebras belongling to this subcategory allows deformations and for this it is necessary to look for a suitable deformation cohohomology which should be constructed using the extra structure. Let \(H^{\ast}_{\mathcal{A}}\) be such a cohomology theory for \(\mathcal{A}\). It should be imminent that the cochain complex defining   \(H^{\ast}_{\mathcal{A}}\) would have to be different from the Hochschild cochain complex. By this we mean that if the \(n\) th cochain space of \(A\in \mathcal{A}\) is \(C^{n}(A,A)= \hom_{k}(A^{\otimes n},A),\) then \(\delta_{\mathcal{A}}\) is not a scalar multiple of \(\delta_{Hoch},\) where \(\delta_{Hoch}\) and \(\delta_{\mathcal{A}}\) are respective coboundary maps from \(C^{n}(A,A)\) to \(C^{n+1}(A,A)\) in Hochschild cochain complex and cochain complex associated to \(H^{\ast}_{\mathcal{A}}.\) The main theme of this note is to illustrate the above question by considering the subcategory \(\mathcal{A}\) of \(Ass\) to be the subcategory of complex Clifford algebras.\newline
			
			Clifford algebras were invented by William K. Clifford who introduced a new multiplication rule into Grassmann's exterior algebra \(\bigwedge \mathbb{R}^{n}.\) A Clifford algebra is a unital associative algebra and generalizes the real numbers, complex numbers and Hamilton's quaternions and plays important roles in geometry and theoretical physics. Clifford algebras can be seen as deformations of exterior algebras. Now it is well known that complex Clifford algebras are rigid i.e. any deformation is equivalent to null deformation with respect to the Hochschild cohomology. The Hochschild cohomology groups of complex Clifford algebras are easy to compute due to the fact that a complex Clifford algebra is either isomorphic to a complex central simple algebra or isomorphic to a direct sum of two isomorphic complex central simple algebras.\newline
			
			In this article we construct a cochain complex of a complex Clifford algebra \(A\) with coefficients in \(A\), for every arbitrary choice of an ordered orthogonal basis \(B\), in a combinatorial fashion which is not of Hochschild cochain complex type. We call the corresponding cohomology by {\it Clifford cohomology} and will be denoted by \(H^{\ast}_{Cl(B)}(A)\). It turns out that deformations of \(A\in \mathcal{A}\), where \(\mathcal{A}\) is the subcategory of complex Clifford algebras of \(Ass,\) is controlled by \(H^{\ast}_{Cl(B)}(A)\) and it also classifies \(\mathcal{A}\) upto Morita equivalence.\newline
			In section \(2\) and \(3\) we briefly recall some facts about Clifford algebras and deformations of associative algebras respectively. In section \(4\) we introduce {\it Clifford cohomology.} In the final section \(5\), we obtained some interesting observations about formal deformations of the algebra of smooth sections of a complex Clifford algebra bundle over an even dimensional orientable Riemannian manifold \(M\) which admits a \(Spin^{c}\) structure by computing its Hochschild cohomology groups.

			\section{Clifford algebra}
			
			We start by recalling some basic facts about Clifford algebras \cite{Lounesto} \cite{Lundholm}. There are lots of literatures on Clifford algebras to be mentioned. Let \(V\) be an \(n\) dimensional complex vector space equipped with a non degenerated quadratic form \(q\) and \(T(V)= \bigoplus _{k= 0}^{\infty} \otimes^{k}V \) be the tensor algebra over \(V.\)
			
			\begin{Def}
				A {\it Clifford algebra} \(C(n)\) over \(V\) is defined to be the quotient algebra \(T(V)/I_{q},\) where \(I_{q}\) is the two-sided ideal generated by elements of the form \(v\otimes v- q(v)\) for all vectors \(v\) in \(V.\) The product in \(C(n)\) is called the {\it Clifford product.}

			\end{Def}
		Now the corresponding bilinear form of \(q\) is \(\beta _{q}(u,v)= \frac{1}{2} (q(u+v)-q(u)-q(v)).\) We note that in \(C(n),\) \[q(u+v)= (u+v)^{2}= u^{2}+ uv+ vu+ v^{2}= q(u)+ uv+ vu+ q(v)\]
		and so \(uv+ vu= 2\beta_{q}(u,v).\) Thus if \(u\) and \(v\) are orthogonal vectors then \(uv= -vu\) in \(C(n).\)\newline
		
		Now if we choose an orthogonal basis \(B=\{v_{1},\cdots, v_{n}\}\) of \(V\) then the set \(\mathcal F\) consisting of \(2^{n}\) many elements given by \(\mathcal F=\{1, E_{m_{1}m_{2}\cdots m_{k}}\}\), \(1\leq k\leq n\) and \(1\leq m_{1}< m_{2}< \cdots <m_{k}\leq n\) form a basis of the vector space \(C(n)\) where \(E_{m_{1}m_{2}\cdots m_{k}}= v_{m_{1}}v_{m_{2}}\cdots v_{m_{k}}\) (here right side is the Clifford product of \(v_{m_{1}},\cdots ,v_{m_{k}}\)) and thus \(\dim C(n)= 2^{n}.\)\newline
		
		\textit{Periodicity of Clifford algebras:}  Let \(M(k)\) denotes the matrix algebra of \(k\times k\) complex matrices. It is known that there is an isomorphism between \(C(n+2)\) and \(C(n)\otimes _{\mathbb{C}} M(2).\) Now from this isomorphism and along with the fact that \(C(0)\cong \mathbb{C}\) and \(C(1)\cong \mathbb{C}\oplus \mathbb{C}\) it follows that:
		\[ 
		\begin{array}{c l}
			C(2n) \cong M(2^{n})\\
			C(2n+1)\cong M(2^{n})\oplus M(2^{n}).\\
		\end{array}
		 \]
		 In paricular, \(C(n+2)\) and \(C(n)\) are Morita equivalent.

		 \section{Deformation of associative algebras and Hochschild cohomology}
		 
		 Let us start with a short review of deformation theory of associative algebras \cite{Markl} \cite{gerstenhaber} \cite{gerstenhaber 2}. Let \(A\) be an associative algebra over a field \(k\) of characteristic zero, \(k[[t]]\) be the formal power series ring over \(k\) and \(A[[t]]\) is the formal power series over \(A\) which is a \(k[[t]]\) algebra.
		 
		 \begin{Def}
		 	A {\it formal deformation}  of \(A\) with base $k[[t]]$ is a \(k[[t]]\)-bilinear multiplication law \(\mu_t :A[[t]]\otimes_{k[[t]]} A[[t]]\rightarrow A[[t]]\) on the spaces \(A[[t]]\)  of formal power series in a variable $t$ with coefficients in $A,$ satisfying the follwoing properties:
 $$\mu_t (a,b)= \mu_{0}(a,b)+ \mu_{1}(a,b)t+ \mu_{2}(a,b)t^{2}+ \cdots \cdots ~~\mbox{for all}~~ a, b\in A,$$
where   \(\mu_{0}(a,b)= ab\)  is the original multiplication in \(A,\) and $\mu_t$ is associative, which is equivalent to the equation 
$$\mu_t(\mu_t (a, b), c) = \mu_t(a, \mu_t(b,c))~~~~\mbox{for}~~a, b, c, \in A,$$
or, equivalently,
		 	
		 	\[\sum_{i+j= k. i,j\geq 0}(\mu_{i}(\mu_{j}(a,b)),c)- \mu_{i}(a,\mu_{j}(b,c))= 0 \]
		 	
		 	for all \(a,\) \(b,\) \(c\in A\) and for each \(k\geq 1.\)\newline
		 	
		 \end{Def}
		 	
		 	If one chooses \(\mu_{i}= 0\) for all \(i\geq 1\) then the deformation of \(A\) is called null deformation.\newline
		 	
		 	
		 	Let \(\mu_t = \mu_{0} + \mu_{1}t+ \mu_{2}t^{2}+ \cdots\) and  \(\mu_t ^\prime= \mu^\prime_{0} + \mu^\prime_{1}t+ \mu^\prime_{2}t^{2}+ \cdots\) be two deformations of \(A.\) Now we say \(\mu_t\) and \(\mu^\prime_t\) are are equivalent if there exists a \(k[[t]]\) linear map \(u : A[[t]]\rightarrow A[[t]]\) defined by \(u= id_{A}+ \phi_{1}t+ \phi_{2}t^{2}+ \cdots ,\) \(\phi_{i}\in Hom_{k} (A,A),\) \(i\in \mathbb{N},\) such that
		 	
		 	\[u\circ \mu_t= \mu^\prime \circ (u\otimes u).\]
		 	
		 	If every formal deformation of \(A\) is equivalent to null deformation then \(A\) is called {\it rigid}.\newline
		 	
		 	The main tool in studing deformation theory of an associative algebra \(A\) is the Hochschild cochain complex \(C^{\ast}(A,A):\)
		 	
		 	\[0\rightarrow C^{0}(A,A)\stackrel{\delta_{Hoch}}{\rightarrow} \cdots \stackrel{\delta_{Hoch}}{\rightarrow}
		 	C^{n}(A,A)\stackrel{\delta_{Hoch}}{\rightarrow} C^{n+1}(A,A) \stackrel{\delta_{Hoch}}{\rightarrow} \cdots \]
		 	
		 	where \(C^{0}(A,A)= A\) and \(C^{n}(A,A)= \hom_{k}(A^{\otimes n},A)\) is the space of Hochschild \(n\)-cochains, i.e., the \(n\)-linear maps \(f\) on a \(A\) with values in \(A.\) The differential \(\delta_{Hoch}: C^{n}(A,A)\rightarrow C^{n+ 1}(A,A)\) is defined by:\vspace{6mm}
		 	
		\[	(\delta_{Hoch} f) (a_{0}\otimes a_{1}\otimes \cdots \otimes a_{n})= a_{0} f(a_{1}\otimes \cdots \otimes a_{n})\]
		\[+ \sum_{i= 1}^{n} (-1)^{i} f(a_{0}\otimes \cdots \otimes a_{i-1}a_{i}\otimes \cdots \otimes a_{n})+ (-1)^{n+ 1} f(a_{0}\otimes \cdots \otimes a_{n- 1})a_{n}.\]
		
		It turns out that \(\delta_{Hoch}^{2}= 0\) and Hochschild cohomology of \(A\) with coefficients in \(A\) is defined by \(HH^{\ast}(A)= H^{\ast} (C^{\ast}(A,A); \delta_{Hoch}).\) Hochschild cohomology groups are invariant under Morita equivalence [for details regarding Hochschild cohomology see \cite{Kassel} \cite{Loday}]. 

It turns out that for a formal deformation of $A$ as defined above the coefficient $\mu_1$ is a Hochschild $2$-cocycle, that is, $\delta_{Hoch} (\mu_1) = 0,$ and is called the infinitesimal of the deformation. 
\begin{Def}
An infinitesimal deformation of \(A\) is a deformation of the form \(A[[t]]/(t^{2})\). More generally, a one parameter deformation of order $n,$  is a deformation with base $k[[t]]/(t^{n+1}),$ given by  $\mu_t$ modulo $t^{n+1}.$  In this case, the associativity condition in the definition above holds for the $2$-cochains $\mu_i$ for $0\leq i \leq n.$ 
\end{Def}

Next comes the question of extending an infinitesimal deformation to a full-blown deformation. If we start with an arbitrary Hochschild $2$-cocycle $\mu_1,$ it need not be an infitesimal of a formal deformation. If it be so, then we say $\mu_1$ is integrable. The integrability of $\mu_1$ implies an infinite sequence of relations which may be interpreted as the vanishing of the obstructions to the integration of $\mu_1.$

Suppose we have a deformation of $A$ of order $n \geq 1$ given by multiplication $\mu_t$ modulo $t^{n+1}$ and we would like to extend this to a deformation of order $n+1.$ Then, $\mu_t$ modulo $t^{n+2}$ must be associative. This gives rise to a $3$-cochain
$$ G (a, b, c) = \Sigma_{i+j=n + 1}\mu_i(\mu_j(a, b), c) - \mu_i(a, \mu_j(b, c)), ~~i>0~~j>0~~\mbox{and}~~a, b, c \in A .$$
It turns out that $G$ is a $3$-cocycle and is called the obstruction cocycle. 
		
We wish to end this section by quoting the following well-known theorems:
		
		\begin{Thm} 
			There is a one-to-one correspondence between the space of equivalence classes of infinitesimal deformations of \(A\) and the second Hochschild cohomology \(HH^{2}(A)\) of \(A\) with coefficients in itself.
		\end{Thm}
	
	\begin{Thm}
		
	Let \(A\) be an associative algebra such that \(HH^{2}(A)= 0.\) Then all formal deformations of \(A\) are equivalent to null deformation, in other words, \(A\) is rigid.
	\end{Thm}
\begin{Thm}
A deformation of $A$ of order $n$ extends to a deformation of order $n+1$ if and only if the cohomology class of the associated obstruction cocycle $G$ vanishes. Thus, if $HH^3(A) = 0$ then, any Hochschild $2$-cocycle is integrable.
\end{Thm}

		\section{Clifford cohomology}
		
		In this section we will construct a cochain complex of a Clifford algebra \(C(n)\) over an \(n\) dimensional complex vector space \(V\) equipped with a non degenerated quadratic form \(q\), which we call \textit{ Clifford cochain complex.} Firstly for every choice of an ordered orthogonal basis \(B\) of \(V\) we will define a bilinear product \(\overline{\circ}_{B}\) on \(C(n)\) and then use it to define the {\it Clifford cochain complex}. It turns out that if \(B\) consists of orthonormal vectors then the {\it Clifford cochain complex} of \(C(n)\) coincides with the Hochschild cochain complex.\newline
		We start with a finite set \(X= \{1,2,\cdots ,n\}.\) Let \(\mathbb{S}\) be the set of all finite sequences in \(X\) and we define a \(``2\)-shuffle", denoted by \(sh_{i},\) on \newline \(\{m_{1},\cdots ,m_{i}, m_{i+1}, \cdots ,m_{k}\}\in \mathbb{S}\)  by
		
		\[\{m_{1},\cdots ,m_{i}, m_{i+1}, \cdots ,m_{k}\}\stackrel{sh_{i}}{\longrightarrow} \{m_{1},\cdots ,m_{i+1}, m_{i}, \cdots ,m_{k}\} \]
		
		for \(1\leq i \leq k-1.\)\newline
		 Now let \(P(X)\uparrow\) be the set of all ordered subsets \(\{m_{1},m_{2},\cdots , m_{k}\}\) of \(X\) such that \(m_{1}< m_{2}< \cdots < m_{k},\) \(1\leq m_{i}\leq n\) for \(i= 1,2,\cdots ,k.\) We consider the empty set, \(\{\quad \},\) as a member of \(P(X)\uparrow .\) Let us define a binary operation \(\circ\) on \(P(X)\uparrow\) as follows:
		
		\[\{m_{1},m_{2},\cdots ,m_{k}\}\circ \{t_{1},t_{2}, \cdots ,t_{s}\}= \{l_{1},l_{2},\cdots ,l_{p}\},\]
		
		if \(\{m_{1},m_{2},\cdots ,m_{k}\}\neq \{t_{1},t_{2}, \cdots ,t_{s}\},\) where \(\{l_{1},l_{2},\cdots ,l_{p}\}\) is obtained by applying the minimum number of \(2\)-shuffles \(\tau\) on the sequence\newline \(\{m_{1},m_{2},\cdots ,m_{k}, t_{1},t_{2},\cdots ,t_{s}\}\) to get a monotonic increasing sequence and then deleting \(m_{i}\) and \(t_{j}\) if \(m_{i}= t_{j},\) \(1\leq i\leq k,\) \(1\leq j\leq s.\)  We note that \(p= card(\{m_{1},m_{2},\cdots ,m_{k}\}\bigtriangleup \{t_{1},t_{2}, \cdots ,t_{s}\})\) (here \(\bigtriangleup\) means symmetric difference of two sets). Let us explain this by an example:\newline
		We take \(X= \{1,2,3,4\}\) and we compute \(\{1,3,4\}\circ \{2,3\}.\) We see \(\{1,3,4\}\circ \{2,3\}= \{1,2,4\}.\) This is obtained by successively applying \(2\)-shuffles \(sh_{3}, sh_{2}\) and \(sh_{4}\) on the sequence \(\{1,3,4,2,3\}\) to get \(\{1,2,3,3,4\}\) and then finally deleting \(3\) and \(3\) to get \(\{1,2,4\}.\) Here \(\tau = 3.\)\newline
		Clearly \(\{m_{1},m_{2},\cdots ,m_{k}\}\circ \{m_{1},m_{2},\cdots ,m_{k}\}= \{\quad \},\) \(1\leq k\leq n,\) is obtained by applying \(\frac{k(k-1)}{2}\) many \(2\)-shuffles on the sequence \(\{m_{1},\cdots ,m_{k},m_{1},\cdots ,m_{k}\}\) and then by deleting everything.\newline
		Let \(B=\{v_{1},\cdots, v_{n}\}\) be an ordered orthogonal basis of an \(n\) dimensional complex vector space \(V\) and the set \(\mathcal F=\{1, E_{m_{1}m_{2}\cdots m_{k}}\}\), where \(E_{m_{1}m_{2}\cdots m_{k}}= v_{m_{1}}v_{m_{2}}\cdots v_{m_{k}}\), \(1\leq k\leq n\), \(1\leq m_{1}< m_{2}< \cdots <m_{k}\leq n\), form a basis of \(C(n).\) Now we define a map \(\circ_{B}: \mathcal F \times \mathcal F \rightarrow C(n)\) by
		\[E_{m_{1}m_{2}\cdots m_{k}} \circ_{B} E_{t_{1}t_{2}\cdots t_{s}}= (-1)^{\tau} E_{l_{1}l_{2}\cdots l_{p}},\] where \(\{l_{1},l_{2},\cdots ,l_{p}\}= \{m_{1},m_{2},\cdots ,m_{k}\}\circ \{t_{1},t_{2},\cdots ,t_{s}\}\) 
		if \(\{m_{1},m_{2},\cdots ,m_{k}\}\neq \{t_{1},t_{2}, \cdots ,t_{s}\},\)
		\[E_{m_{1}m_{2}\cdots m_{k}}\circ_{B} E_{m_{1}m_{2}\cdots m_{k}}= (-1)^{\frac{k(k-1)}{2}}\]
		and
		\[ 1\circ_{B} E_{m_{1}m_{2}\cdots m_{k}}= E_{m_{1}m_{2}\cdots m_{k}}= E_{m_{1}m_{2}\cdots m_{k}}\circ_{B} 1.\] 
		Finally we extend \(\circ_{B}\) bilinearly on \(C(n)\times C(n)\) and get a linear map \(\overline{\circ}_{B}: C(n)\otimes C(n)\rightarrow C(n).\)
		
		\begin{Lem}
			If \(B'= \{e_{1},e_{2},\cdots ,e_{n}\}\) is an ordered orthonormal basis of \(V\) then \(\overline{\circ}_{B'}\) coincides with the Clifford product in \(C(n).\)
		\end{Lem}
	
	\begin{proof}
		Let \(\overline{E}_{m_{1}m_{2}\cdots m_{k}}= e_{m_{1}}e_{m_{2}}\cdots e_{m_{k}}.\) Now the lemma follows from the fact that
		\[\overline{E}_{m_{1}m_{2}\cdots m_{k}} \circ_{B'} \overline{E}_{t_{1}t_{2}\cdots t_{s}}= (-1)^{\tau} \overline{E}_{l_{1}l_{2}\cdots l_{p}}=  e_{m_{1}}e_{m_{2}}\cdots e_{m_{k}} e_{t_{1}}e_{t_{2}}\cdots e_{t_{s}}. \]
		if \(\{m_{1},\cdots ,m_{k}\}\neq \{t_{1}, \cdots ,t_{s}\},\)
		\[\overline{E}_{m_{1}m_{2}\cdots m_{k}}\circ_{B'} \overline{E}_{m_{1}m_{2}\cdots m_{k}}= (-1)^{\frac{k(k-1)}{2}}= e_{m_{1}}e_{m_{2}}\cdots e_{m_{k}}e_{m_{1}}e_{m_{2}}\cdots e_{m_{k}} \]
		and
		\[ 1\circ_{B'} \overline{E}_{m_{1}m_{2}\cdots m_{k}}= \overline{E}_{m_{1}m_{2}\cdots m_{k}}= e_{m_{1}}e_{m_{2}}\cdots e_{m_{k}}= \overline{E}_{m_{1}m_{2}\cdots m_{k}}\circ_{B'} 1.\]

	\end{proof}

Now let \(C^{q}(C(n))= Hom_{\mathbb{C}}(C(n)^{\otimes q},C(n)),\) for all \(q\geq 0\), \(C^{0}(C(n))= C(n)\) and \(\delta_{Cl(B)}: C^{q}(C(n))\rightarrow C^{q+1}(C(n)) \) be the coboundary map defined by
\[	(\delta_{Cl(B)} f) (a_{0}\otimes a_{1}\otimes \cdots \otimes a_{q})= a_{0}\overline{\circ}_{B} f(a_{1}\otimes \cdots \otimes a_{q})\]
\[+ \sum_{i= 1}^{q} (-1)^{i} f(a_{0}\otimes \cdots \otimes a_{i-1}\overline{\circ}_{B} a_{i}\otimes \cdots \otimes a_{q})+ (-1)^{q+ 1} f(a_{0}\otimes \cdots \otimes a_{q- 1})\overline{\circ}_{B} a_{q},\]
where \(f\in C^{q}(C(n)).\) We note that if \(B'= \{e_{1},e_{2},\cdots ,e_{n}\}\) is an ordered orthonormal basis of \(V\) then \(\delta_{Cl(B')}= \delta_{Hoch}\) by Lemma 4.1.

\begin{Lem}
For any choice of an ordered orthogonal basis \(B\) of \(V,\) \(\delta_{Cl(B)}^{2}= 0.\)
\end{Lem}

\begin{proof}
	Let \(B=\{v_{1},\cdots, v_{n}\}\) be an ordered orthogonal basis and \(B'= \{e_{1},e_{2},\cdots ,e_{n}\}\) is an ordered orthonormal basis of \(V\), \(E_{m_{1}m_{2}\cdots m_{k}}= v_{m_{1}}v_{m_{2}}\) \(\cdots v_{m_{k}}\) and \(\overline{E}_{m_{1}m_{2}\cdots m_{k}}= e_{m_{1}}e_{m_{2}}\cdots e_{m_{k}}.\) We take the vector space isomorphism \(\psi :C(n)\rightarrow C(n)\) defined by
	\[\psi (\overline{E}_{m_{1}m_{2}\cdots m_{k}})= E_{m_{1}m_{2}\cdots m_{k}}\]
	
	\[\psi (1)= 1\]
	and consider the isomorphism \(\psi^{\ast} :C^{q}(C(n))\rightarrow C^{q}(C(n))\) defined by
	\[f\longmapsto \psi^{\ast} f,\]
	where
	\[\psi^{\ast}f(a_{0}\otimes a_{1}\otimes \cdots \otimes a_{q-1})= \psi^{-1}[f(\psi^{\otimes q}(a_{0}\otimes a_{1}\otimes \cdots \otimes a_{q-1}))].\]
	It is noted earlier that \(\delta_{Cl(B')}= \delta_{Hoch}\) and by Lemma 4.1, \(\overline{\circ}_{B'}\) coincides with the Clifford product in \(C(n).\) Now it is an easy check that the following diagram


$$\xymatrix{C^{q}(C(n)) \ar[r]^{\delta_{Cl(B)}}  \ar[d]^{\psi^{\ast}}
		& C^{q+1}(C(n)) \ar[d]^{\psi^{\ast}} \\
		C^{q}(C(n)) \ar[r]^{\delta_{Hoch}}
		& C^{q+1}(C(n))}$$	
	
commutes for all \(q\geq 0\) and therefore \(\psi^{\ast}(\delta_{Cl(B)}^{2}f)= \delta_{Hoch}^{2} (\psi^{\ast}f)= 0,\) \(f\in C^{q}(C(n)).\) Finally as \(\psi^{\ast}\) is an isomorphism so \(\delta_{Cl(B)}^{2}f= 0.\) This completes the proof.
	
\end{proof}

We define the {\it Clifford cochain complex} \((C^{\ast}(C(n)); \delta_{Cl(B)})\) of the Clifford algebra \(C(n)\) associated to an ordered orthogonal base \(B\) of \(V\) by
\[0\rightarrow C(n)\stackrel{\delta_{Cl(B)}}{\rightarrow}\cdots \stackrel{\delta_{Cl(B)}}{\rightarrow} C^{q}(C(n))\stackrel{\delta_{Cl(B)}}{\rightarrow} C^{q+1}(C(n))\stackrel{\delta_{Cl(B)}}{\rightarrow} \cdots \]
and the {\it Clifford cohomology} of \(C(n)\) associated to \(B\) by
\[H^{\ast}_{Cl(B)}(C(n))= H^{\ast}(C^{\ast}(C(n)); \delta_{Cl(B)}).\]

\begin{Thm}
	
	For any ordered orthogonal base \(B\) of \(V,\) \(H^{\ast}_{Cl(B)}(C(n))\cong HH^{\ast}(C(n)).\)
	
\end{Thm}

\begin{proof}
	The proof at once follows from the commutativity of the diagram in Lemma 4.2 along with the fact that \(\psi^{\ast}\) is an isomorphism.
\end{proof}

\begin{Cor}
	
Up to isomorphism Clifford cohomology groups of \(C(n)\) are independent of the choice of an ordered orthogonal base of \(V.\)
\end{Cor}
\begin{proof}
	It readily follows from Theorem 4.3.
\end{proof}

\begin{Rem}
	
	It follows from the construction of {\it Clifford cochain complex} that in general the coboundary maps \(\delta_{Cl(B)}\) are not scalar multiples of \(\delta_{Hoch}\) for arbitrary choices of an ordered orthogonal basis \(B\) and consequently it is not of Hochschild cochain complex type while {\it Clifford cohomology} being isomorphic to Hochschild cohomology, controls the deformations of complex Clifford algebras.
\end{Rem}

Let \(B= \{v_{1},\cdots ,v_{n}\}\) is an ordered orthigonal basis of an \(n\) dimensional complex vector space \(V.\)

\begin{Prop}
	If \(n\) is odd, then \(v_{1}v_{2}\cdots v_{n}\) is a \(0\) cocycle in the Clifford cochain complex associated to \(B.\)

\end{Prop}

\begin{proof}
	If \(m\in C(n)\) then \(\delta_{Cl(B)}(m)(a)= m \overline{\circ}_{B} a- a\overline{\circ}_{B} m\) for all \(a\in C(n).\) Now if \(n\) is odd then we note that \(v_{1}v_{2}\cdots v_{n} \overline{\circ}_{B} a = a \overline{\circ}_{B} v_{1}v_{2}\cdots v_{n}\) for all \(a\in C(n)\) and consequently \(v_{1}v_{2}\cdots v_{n}\) is a \(0\) cocycle in the Clifford cochain complex associated to \(B.\)
\end{proof}

\begin{Prop}
	\begin{enumerate}
	If \(n\) is odd then
	
	\[ H^{i}_{Cl(B)}(C(n)) = \left\lbrace
	\begin{array}{c l}
		\mathbb{C}\oplus \mathbb{C} & \text{if \(i= 0\) }\\
		0 & \text{if \(i> 0\)}\\
	\end{array}
	\right. \]
	and if \(n\) is even then
	
		\[ H^{i}_{Cl(B)}(C(n)) = \left\lbrace
	\begin{array}{c l}
		\mathbb{C} & \text{if \(i= 0\) }\\
		0 & \text{if \(i> 0\)}\\
	\end{array}
	\right. \]
\end{enumerate}
\end{Prop}

\begin{proof}
	First we note that \(H^{\ast}(C(n))\cong HH^{\ast}(C(n))\) [Theorem 4.3]. Let \(n\) be odd. Then \(C(n)\) is Morita equivalent to \(C(1)\cong \mathbb{C}\oplus \mathbb{C}.\) As Hochschild cohomology is invariant under Morita equivalence therefore \(HH^{\ast}(C(n))\cong HH^{\ast}(\mathbb{C}\oplus \mathbb{C})\cong HH^{\ast}(\mathbb{C})\oplus HH^{\ast}(\mathbb{C}).\) Now if \(n\) is even then \(C(n)\) is Morita equivalent to \(C(0)\cong \mathbb{C}\) and \(HH^{\ast}(C(n))\cong HH^{\ast}(\mathbb{C}).\) Again it is known that for any field \(k,\) \(HH^{0}(k)= k\) and \(HH^{i}(k)= 0\) for \(i> 0\) (\cite{Loday}, 1.5.5) and this completes the proof.

\end{proof}

We end this section by showing that in the category of complex Clifford algebras, {\it Clifford cohomology} associated to any ordered orthogonal basis can classify algebras up to Morita equivalence.

\begin{Thm}
Let \(W,\) \(V\) are complex vector spaces of dimension \(m\) and \(n\) with \(B,\) \(B'\) are any two ordered orthogonal basis of them respectively. Then the Clifford algebras \(C(m)\) and \(C(n)\) over \(W\) and \(V\) are Morita equivalent if and only if \(H^{\ast}_{Cl(B)}(C(m))\cong H^{\ast}_{Cl(B')}(C(n)).\)
\end{Thm}

\begin{proof}
	If \(C(m)\) and \(C(n)\) are Morita equivalent then \(HH^{\ast}(C(m))\cong HH^{\ast}(C(n))\) and it follows from Theorem 4.3 that \(H^{\ast}_{Cl(B)}(C(m))\cong H^{\ast}_{Cl(B')}(C(n)).\) \newline
	Conversely, let \(H^{\ast}_{Cl(B)}(C(m))\cong H^{\ast}_{Cl(B')}(C(n)).\) Now it follows from Propsition 4.7 that \(m\) and \(n\) must be both even or both odd and Consequently \(C(m)\) and \(C(n)\) are Morita equivalent.
\end{proof}

\section{Formal deformations of smooth sections of complex Clifford algebra bundle}
The aim of this last section is to study  Hochschild cohomology groups and formal deformations of the algebra of smooth sections of a complex Clifford algebra bundle over an even dimensional orientable Riemannian manifold \(M\) which admits a \(Spin^{c}\) structure. It turns out that if \(M\) is \(2\)-dimensional then the algebra of smooth sections of the complex Clifford algebra bundle over \(M\) is highly non-rigid in the sense that it admits infinitely many inequivalent formal deformations.\newline
			
			We start by recalling very briefly some facts about \(Spin\) and \(Spin^{c}\) manifolds (for details see \cite{Lawson}). 
\begin{Def}
Let \(M\) be an orientable Riemannian manifold of dimension \(n\), \(P_{SO}(M)\) be the oriented orthonormal frame bundle over \(M\) and \(Spin(n)\) is the double covering group of \(So(n)\). The manifold \(M\) is said to have a \(Spin\) structure if there exists a \(Spin(n)\) bundle \(P_{Spin}(M)\) over \(M\) and an equivariant bundle map: \(P_{Spin}(M)\rightarrow P_{SO}(M)\). 
\end{Def}
The complex analogue of \(Spin(n)\) group is \(Spin^{c}(n)= Spin(n)\times_{\mathbb{Z}_{2}} U(1) \subset R(n)\otimes \mathbb{C}\), where \(R(n)\) denotes the real Clifford algebra over \(\mathbb{R}^{n}\) equipped with a positive definite form. 
\begin{Def}
We say  \(M\) admits a \(Spin^{c}\) structure if there exists a \(Spin^{c}\) bundle \(P_{Spin^{c}}(M)\) over \(M\), a \(U(1)\) bundle \(P_{U}(M)\) over \(M\) and an equivariant bundle map: \(P_{Spin^{c}}(M)\rightarrow P_{SO}(M)\times P_{U}(M)\). 
\end{Def}
\begin{Def}
Let $ k = \mathbb R ~~\mbox{or}~~ \mathbb C,$ and $E\rightarrow M$ be a smooth $k$- vector bundle over a manifold $M.$ Then $E$ is called a bundle of $k$-algebras, if each fibre $E_x$ is a $k$-algebra for any $x \in M,$ such that the algebra operations are smooth.
\end{Def}

\begin{Def}
Suppose $E\rightarrow M$ be a smooth $k$- vector bundle over a manifold $M.$ A $k$-vector bundle \textbf{S} over \(M\) is said to be  bundle of  $E$-module if there is a smooth bundle map: \(E\otimes \textbf{S}\rightarrow \textbf{S}\) that makes \(\textbf{S}_x\) a $E_x$-module for each $x \in M.$
\end{Def}

\begin{Def}
The Clifford algebra bundle \(Cl(M)\) over a smooth manifold \(M\) is obtained from the tangent bundle $TM$ by replacing each fibre \(T_xM\) by the Clifford algebra over \(T_xM\). More precisely,  the  \(Cl(M)= \cup_{x\in M}Cl(T_xM)\). The complex Clifford algebra bundle \(\mathbb{C}l(M)\) is obtained from $Cl(M)$  by complexifying each fibre, that is, \(\mathbb{C}l(M)= Cl(M)\otimes \mathbb{C}\).
\end{Def}

\begin{Rem}
It is well-known that the real Clifford algebra bundle \(Cl(M)\) over \(M\) is \(P_{SO}(M)\times_{So(n)} R(n)\) where \(So(n)\rightarrow Aut(R(n))\) is the natural action. Moreover, if \(M\) has a \(Spin\) structure then \(Cl(M)\) can also be expressed as \(P_{Spin}(M)\times_{Spin(n)} R(n)\) where the action of \(Spin(n)\) on \(R(n)\) is the adjoint action.  If \(M\) is \(Spin^{c}\) then \(\mathbb{C}l(M)\) can also be obtained as \(P_{Spin}(M)\times_{Spin^{c}(n)}R(n)\otimes \mathbb{C}\), where the action of \(Spin^{c}(n)\) on \(R(n)\otimes \mathbb{C}\) is the adjoint action.
\end{Rem}

			
			From now on we denote the algebra of smooth complex functions on \(M\) by \(C^{\infty}(M)\) and the algebra of smooth complex sections of a complex vector bundle \(E\) over\(M\) by \(\Gamma ^{\infty}(E)\). Clearly \(\Gamma^{\infty} (\mathbb{C}l(M))\) is an associatiove unital algebra over \(\mathbb{C}\).
			
			\begin{Thm}
				If \(M\) is an orientable Riemannian manifold of dimension \(2m\), \(m\geqq 1\), which admits a \(Spin^{c}\) structure then the Hochschild cohomology group \(HH^{k}(\Gamma^{\infty} (\mathbb{C}l(M)))\) is non trivial for \(k\leq 2m\) and trivial for \(k> 2m\).
			\end{Thm}
			
			\begin{proof}
				First we note that as \(M\) is of dimension \(2m\), therefore the complex Clifford algebra bundle \(\mathbb{C}l(M)\) is a bundle of complex matrix algebra \(M(2^{m})\) over \(M\). Moreover the existence of a \(Spin^{c}\) structure on \(M\) ensures that there is a complex vector bundle \(E\) of fiber dimension \(2^{m}\) over \(M\) which is a \(\mathbb{C}l(M)\) module i.e. there is a continuous bundle map: \(\mathbb{C}l(M)\otimes E\rightarrow E\) (see \cite{Lawson}, Proposition II.3.8 for real version).\newline
				Now by the smooth version of Serre-Swan's theorem (\cite{Nes}, Theorem 11.32), \(\Gamma ^{\infty}(E)\) is a  finitely generated and projective module over \(C^{\infty}(M)\). Also by Morita's theorem (see e.g. \cite{Lam}, Sec.18; \cite{Burs}, Theorem 4.1.) it follows that \(\Gamma^{\infty}(End (E))\) (\(End(E)\) is the endomorphism bundle) is Morita equivalent to \(C^{\infty}(M)\) where \(\Gamma^{\infty}(E)\) is an invertible (\(\Gamma^{\infty}(End(E))\), \(C^{\infty}(M)\))- bimodule (see \cite{Burs}, Example 4.2.). Again \(E\) being a \(\mathbb{C}l(M)\) module and as \(\mathbb{C}l(M)\) is a bundle of complex matrix algebra \(M(2^{m})\), therefore \(End(E)\cong \mathbb{C}l(M)\) and consequently \(\Gamma ^{\infty}(\mathbb{C}l(M))\) is Morita equivalent to \(C^{\infty}(M)\).\newline
				It is known that Morita equivalent algebras have isomorphic Hochschild cohomology groups. Let \(\bigwedge ^{k}TM\otimes \mathbb{C}\) denotes the complexified \(k\)-th exterior budle over \(M\) and let \(\chi ^{\ast}(M)= \bigoplus _{k= 0}^{\infty}\Gamma ^{\infty}(\bigwedge ^{k}TM\otimes \mathbb{C})\).
				Now we consider the Hochschild-Kostant-Rosenberg map \(\textit{U}: \chi^{k}(M)\rightarrow HH^{k}(C^{\infty}(M))\) defined by: \(X\mapsto \textit{U}(X)\), \(X\in \chi^{k}(M)= \Gamma ^{\infty}(\bigwedge ^{k}TM\otimes \mathbb{C})\) and \(\textit{U}(X)(f_{1},\cdots ,f_{k})= \frac{1}{k!} X(df_{1},\cdots ,df_{k})\) [\(X\) can be viewed as a multilinear alternating map: \(X: \Omega^{1}_{c}(M)\times \cdots \times \Omega^{1}_{c}(M)\rightarrow C^{\infty}(M)\), where \(\Omega^{1}_{c}(M)\) is the space of complexified \(1\) forms i.e. \(\Omega^{1}_{c}(M)= \Gamma^{\infty}(T^{\ast}M\otimes \mathbb{C})\), \(T^{\ast}(M)\) is the cotangent bundle of \(M\)]. We note that if \(X= X_{1}\wedge \cdots \wedge X_{i}\wedge \cdots \wedge X_{k}\), where \(X_{i}\in \chi^{1}(M)\), then \(\textit{U}(X)(f_{1},\cdots ,f_{k})= \frac{1}{k!}\sum_{\sigma \in S_{k}}(\mbox{Sign}~\sigma)X_{\sigma (1)}(f_{1})\cdots X_{\sigma (k)}(f_{k})\). \newline 
				Finally as the Hochschild-Kostant-Rosenberg map \(\textit{U}: \chi^{k}(M)\rightarrow HH^{k}(C^{\infty}(M))\) is injective (\cite{Wal}, Cor: 6.2.47) and there are infinitely elements in \(\chi^{k}(M)\) for \(k\leq 2m\), therefore the non triviality of \(HH^{k}(\Gamma^{\infty}(\mathbb{C}l(M)))\), \(k\leq 2m\) follows from this. The triviality of \(HH^{k}(\Gamma^{\infty}(\mathbb{C}l(M)))\) while \(k> 2m\) follows from the fact that \(\bigwedge^{k}T_{p}M\otimes \mathbb{C}\) is the zero vector space for \(k> 2m\) and for all \(p\in M\).
				
			\end{proof}

			A star product on \(M\) is a formal deformation of \(C^{\infty}(M)\), i.e. an associative product \(\star\) on the \(\mathbb{C}[[t]]\) module \(C^{\infty}(M)[[t]]\) given by: for \(f\), \(g\in C^{\infty}(M)\),
			
			\[f\star g= fg+ \sum_{i=1}^{\infty} \mu_{i}(f,g)t^{i}\]
			
			where \(\mu_{i}: \C^{\infty}(M)\times C^{\infty}(M)\rightarrow C^{\infty}(M)\), \(i= 1,2,\cdots\) are bi-differential operators. We denote the equivalence classes of star products by \(Def(M)\).\newline
			
			Any Lie bracket \(\{\cdot ,\cdot\}\) on \(C^{\infty}(M)\) which is compatible with the pointwise product on \(C^{\infty}(M)\) via the Liebniz rule is called a Poisson structure on \(M\). Given any star product \(\star\) on \(M\), it is known that \(\{f,g\}= \frac{1}{i}(\mu_{1}(f,g)- \mu_{1}(g,f))\); \(f\), \(g\in C^{\infty}(M)\), is a Poisson structure on \(M\) (see \cite{Burs} Sec. 3.2.).
			
			\begin{Thm}
				If \(M\) is a \(2\)-dimensional \(Spin^{c}\) manifold then \(\Gamma^{\infty}(\mathbb{C}l(M))\) admits infinitely many inequivalent formal deformations.
			\end{Thm}
			
			\begin{proof}
				As \(M\) is \(Spin^{c}\) therefore \(\Gamma^{\infty}(\mathbb{C}l(M))\) is Morita equivalent to \(C^{\infty}(M)\) (follows from the proof of Theorem 5.7.). Again as the set of equivalence classes of formal deformations is Morita invariant (\cite{gerstenhaber}, section 16) so it suffices to explore \(Def(M)\).\newline
				As \(M\) is of dimension \(2\), therefore each complex bi-vector field  \(\pi \in \chi^{2}(M)\)\newline\(=\Gamma^{\infty}(\bigwedge^{2}TM \otimes \mathbb{C})\) induces a Poisson structure \(\{\cdot ,\cdot\}_{\pi}\) on \(M\) defined by: \(\{f,g\}_{\pi}= \pi (df,dg)\). Now by Kontsevich's classification result (\cite{Kont}; \cite{Burs} Theorem 3.3.) distinct Poisson structures on \(M\) corresponds to distinct elements in \(Def(M)\). Finally as \(\chi^{2}(M)\) is clearly an infinite set therefore \(Def(M)\) is also infinite. This completes the proof.
			\end{proof}
			
			\begin{Cor}
				\(\Gamma^{\infty}(\mathbb{C}l(S^{1}\times S^{1}))\) and \(\Gamma^{\infty}(\mathbb{C}l(S^{2}))\) have infinitely many inequivalent formal deformations.
			\end{Cor}
			
			\begin{proof}
				
				As \(S^{1}\times S^{1}\) is parallelizable so it is \(Spin^{c}\) and tangent bundle of \(S^{2}\) being stably trivial, it is a \(Spin\) and therefore \(Spin^{c}\). Now the proof readily follows from Theorem 5.8.
			\end{proof}

		\textbf{Acknowledgement.} The authors thank Professor Stefan Waldmann for his valuable suggestions and comments.

		\end{document}